\documentclass[sigconf, nonacm]{acmart}
\usepackage{amsmath, pifont} 
\usepackage{algorithm} 
\usepackage{algpseudocode} 
\usepackage{graphicx}
\usepackage{textcomp}
\usepackage{xcolor}
\usepackage{array}
\usepackage{makecell}
\usepackage{threeparttable}
\usepackage{multirow}
\usepackage{balance}

\newcommand\vldbdoi{10.14778/3685800.3685829}
\newcommand\vldbpages{4090 - 4103}
\newcommand\vldbvolume{17}
\newcommand\vldbissue{12}
\newcommand\vldbyear{2024}
\newcommand\vldbauthors{\authors}
\newcommand\vldbtitle{\shorttitle} 
\newcommand\vldbavailabilityurl{}
\newcommand\vldbpagestyle{empty}

\begin{document}

\title{OptScaler: A Collaborative Framework for Robust Autoscaling in the Cloud}

\author{Ding Zou}
\authornote{Equal contribution}
\affiliation{
\institution{Zhejiang University}
\institution{Ant Group}
}
\email{zoud@zju.edu.cn}

\author{Wei Lu}
\authornotemark[1]
\affiliation{
\institution{Ant Group}
}
\email{xiaobo.lw}
\email{@antgroup.com}

\author{Zhibo Zhu}
\affiliation{
\institution{Ant Group}
}
\email{gavin.zzb}
\email{@antgroup.com}

\author{Xingyu Lu}
\authornote{Corresponding authors}
\affiliation{
\institution{Ant Group}
}
\email{sing.lxy}
\email{@antgroup.com}

\author{Jun Zhou}
\authornotemark[2]
\affiliation{
\institution{Ant Group}
}
\email{jun.zhoujun}
\email{@antgroup.com}

\author{Xiaojin Wang}
\affiliation{
\institution{Ant Group}
}
\email{yueying.wxj}
\email{@antgroup.com}

\author{Kangyu Liu}
\affiliation{
\institution{Ant Group}
}
\email{liukangyu.lky}
\email{@antgroup.com}

\author{Kefan Wang}
\affiliation{
\institution{Ant Group}
}
\email{kefan.wkf}
\email{@antgroup.com}

\author{Renen Sun}
\affiliation{
\institution{Ant Group}
}
\email{renen.sun}
\email{@antgroup.com}

\author{Haiqing Wang}
\affiliation{
\institution{Ant Group}
}
\email{wanghaiqing.whq}
\email{@antgroup.com}

\begin{abstract}
Autoscaling is a critical mechanism in cloud computing, enabling the autonomous adjustment of computing resources in response to dynamic workloads. This is particularly valuable for co-located, long-running applications with diverse workload patterns. The primary objective of autoscaling is to regulate resource utilization at a desired level, effectively balancing the need for resource optimization with the fulfillment of Service Level Objectives (SLOs). 
Many existing proactive autoscaling frameworks may encounter prediction deviations arising from the frequent fluctuations of cloud workloads. Reactive frameworks, on the other hand, rely on real-time system feedback, but their hysteretic nature could lead to violations of stringent SLOs. Hybrid frameworks, while prevalent, often feature independently functioning proactive and reactive modules, potentially leading to incompatibility and undermining the overall decision-making efficacy.
In addressing these challenges, we propose OptScaler, a collaborative autoscaling framework that integrates proactive and reactive modules through an optimization module. The proactive module delivers reliable future workload predictions to the optimization module, while the reactive module offers a self-tuning estimator for real-time updates. By embedding a Model Predictive Control (MPC) mechanism and chance constraints into the optimization module, we further enhance its robustness.
Numerical results have demonstrated the superiority of our workload prediction model and the collaborative framework, leading to over a 36\% reduction in SLO violations compared to prevalent reactive, proactive, or hybrid autoscalers. Notably, OptScaler has been successfully deployed at Alipay, providing autoscaling support for the world-leading payment platform.
\end{abstract}

\maketitle

\pagestyle{\vldbpagestyle}
\begingroup\small\noindent\raggedright\textbf{PVLDB Reference Format:}\\
\vldbauthors. \vldbtitle. PVLDB, \vldbvolume(\vldbissue): \vldbpages, \vldbyear.\\
\href{https://doi.org/\vldbdoi}{doi:\vldbdoi}
\endgroup
\begingroup
\renewcommand\thefootnote{}\footnote{\noindent
This work is licensed under the Creative Commons BY-NC-ND 4.0 International License. Visit \url{https://creativecommons.org/licenses/by-nc-nd/4.0/} to view a copy of this license. For any use beyond those covered by this license, obtain permission by emailing \href{mailto:info@vldb.org}{info@vldb.org}. Copyright is held by the owner/author(s). Publication rights licensed to the VLDB Endowment. \\
\raggedright Proceedings of the VLDB Endowment, Vol. \vldbvolume, No. \vldbissue\ %
ISSN 2150-8097. \\
\href{https://doi.org/\vldbdoi}{doi:\vldbdoi} \\
}\addtocounter{footnote}{-1}\endgroup

\ifdefempty{\vldbavailabilityurl}{}{
\vspace{.3cm}
\begingroup\small\noindent\raggedright\textbf{PVLDB Artifact Availability:}\\
The source code, data, and/or other artifacts have been made available at \url{\vldbavailabilityurl}.
\endgroup
}
\section{Introduction} \label{sec: Intro}
The rapid growth of cloud computing has generated a significant demand for computing resources such as CPU cores and memory \cite{straesser2022not}. As a result, efficient resource management and optimization techniques have received increasing attention to strike a balance between the cost of resources and the strict Service Level Objectives (SLOs) in cloud services. Traditionally, to address drastic workload changes and ensure SLO adherence, resources are often provisioned based on peak demand, leading to substantial resource waste \cite{ding2019characterizing}. In response to this challenge, autoscaling has emerged as a fundamental capability of cloud infrastructure, allowing for the automatic and dynamic scaling of resources in reaction to workload fluctuations \cite{chen2018survey}. With autoscaling, cloud providers can swiftly adapt to varying workloads without manual intervention, thereby achieving optimal performance. This capability not only provides substantial cost savings but also enhances user experience. There are two types of scaling based on how resources are adjusted: vertical scaling, which modifies the resource capacity within existing cluster nodes \cite{farokhi2016hybrid}, and horizontal scaling, which involves adding or removing deployed nodes \cite{zhou2023ahpa}. Major cloud vendors widely favor horizontal scaling due to its ease of implementation and its ability to enhance the availability and fault tolerance of applications \cite{chen2018survey,deng2018establishment,deng2019distributed}.

In this paper, we focus on horizontal scaling in the context of co-located long-running applications (LRAs) with diverse workload patterns. LRAs, also known as online services, are integral in many cloud computing scenarios such as online marketing and content recommendation. These scenarios frequently exhibit time-varying workloads due to periodic user arrivals. LRAs differ from batch jobs, where each job is assigned to an exclusive node that is terminated after processing the job \cite{qian2022robustscaler}. Existing autoscaling methods (also referred to as autoscalers) for LRAs can be classified as either reactive or proactive, based on the timing of scaling, and both have been widely researched \cite{al2017elasticity}. Reactive autoscalers, driven by real-time system feedback \cite{49174, qiu2020firm}, adaptively adjust resources at run-time after unexpected outcomes occur and are favored for their simplicity of implementation. However, real-world workloads often exhibit periodic fluctuations \cite{gambi2013testing}, and the reactive autoscaling's hysteretic nature could lead to violations of SLOs.

In this regard, \cite{podolskiy2018iaas} experimented with the autoscaling methods from Amazon, Microsoft, and Google, leading to the conclusion that their reactive methods had performance pitfalls.
In contrast, proactive autoscalers, which anticipate future workload and scale resources beforehand, might avoid the pitfalls of reactive scaling \cite{Dezhabad2018, xue2022meta}. However, proactive autoscalers with sophisticated prediction techniques may struggle to capture the full picture of uncertain LRAs' workloads, resulting in prediction errors and resource wastage. To enhance the performance of proactive autoscaling under uncertainty (e.g., unpredictable bursts or anomalies), integrating proactive and reactive methods to build a hybrid autoscaler is essential, as highlighted in \cite{straesser2022not}. An ideal paradigm involves the proactive method handling foreseeable workload patterns while the reactive method corrects unexpected deviations.

In practical applications, hybrid autoscalers face their own set of challenges. Firstly, the use of proactive and reactive modules in existing hybrid autoscalers often leads to independent operation, potentially resulting in conflicting scaling decisions when the outputs of these modules differ. While an enhanced prediction model in the proactive module could mitigate these conflicts, \cite{straesser2022not} stresses that reconciling  conflicting outputs to reach a final implementation decision remains a significant challenge in production.
Secondly, the intricate nature of the cloud environment poses challenges in developing a robust hybrid autoscaler. For example:
1) Hardware speed limitations can constrain scaling capacity, possibly causing disparities between intended and actual scaling decisions, ultimately undermining autoscaling performance;
2) Systematic noise in scaling metrics, such as CPU utilization, has the potential to mislead both proactive and reactive scaling decisions;
3) The co-location of LRAs frequently leads to resource contention, further exacerbating the uncertainty in scaling metrics.

The primary goal of this study is to develop an autoscaling framework for LRAs that effectively tackles the challenges previously mentioned. To confront the first challenge, we introduce a novel \textbf{collaborative} autoscaling framework, an advancement of existing hybrid frameworks. This includes an innovative proactive module with enhanced prediction capabilities, a real-time reactive module to address scaling errors stemming from prediction inaccuracies, and an optimization module to manage the trade-off between resource costs and SLO satisfaction.
The term \textbf{collaborative} indicates that our framework orchestrates proactive and reactive modules to make holistic final decisions, thus avoiding potential conflicts between the two modules. To achieve this, our collaborative framework constructs the optimization module with dynamic inputs from both proactive and reactive modules.
The optimization objective is to assist the cloud system in attaining the desired scaling metrics under dynamic workloads. The interpretable optimization objective and constraints (i.e., practical restrictions in the cloud) also improve the interpretability of our framework, an aspect often lacking in many black-box or machine-learning-based autoscalers. Enhanced transparency within the framework can further establish user trust and facilitate potential upgrades of autoscalers \cite{straesser2022not}.

To tackle the second challenge, we explicitly model the essential components of the complex cloud system to enable robust scaling decisions. Specifically:
1) We integrate Model Predictive Control (MPC) into our optimization module to address the speed limitation of scaling. In contrast to traditional control methods such as Proportional-Integral-Derivative (PID) \cite{9599238}, MPC is recognized for its superior potency. In our context, MPC enables optimization of the scaling decision at the current time while accounting for potential workload bursts in future time intervals;
2) We mitigate the impact of noises in scaling metrics by employing the chance-constraint method \cite{shapiro2021lectures}, further enhancing the robustness of MPC;
3) We develop an adaptive estimator within the reactive module, continually adjusted based on real-time system feedback, to monitor the collaborative impact of co-located LRAs on scaling metrics.

The enhancement and collaboration of each module in OptScaler has proven highly effective for managing co-located LRAs with diverse workload patterns. To support this, the following questions are thoroughly investigated:

\begin{enumerate}
\item \textbf{Question 1}: Can OptScaler offer superior workload prediction ability compared to prevalent methods?
\item \textbf{Question 2}: Can OptScaler make more robust scaling decisions than mainstream autoscaling frameworks (e.g., reactive, proactive, and hybrid frameworks) across various real-world workload patterns?
\item \textbf{Question 3}: Taking the autoscaling framework as a whole, what is the overall advantage of OptScaler compared to state-of-the-art autoscalers?
\end{enumerate}

\textbf{Contributions.} Addressing the above questions, OptScaler contributes to cloud resource management in three key aspects:

\begin{enumerate}
\item OptScaler develops an innovative proactive module that surpasses state-of-the-art prediction methods when tested on challenging public and internal workload traces;
\item OptScaler creatively leverages an optimization module to enable the proactive module to collaborate with the reactive module in making final scaling decisions. The incorporation of MPC mechanism and chance constraints into the optimization module enhances its robustness;
\item OptScaler demonstrates its superiority as a comprehensive autoscaling framework, reducing over $36\%$ more SLO violations compared to a state-of-the-art hybrid framework. OptScaler has successfully supported the online autoscaling of LRAs at Alipay, the world-leading payment platform.
\end{enumerate}

The remainder of the paper is organized as follows. Section \ref{sec: Literature} reviews related works. Section \ref{sec: Preliminaries} provides background information about our proposed autoscaling framework. Section \ref{sec: Proposed} elaborates on the proposed framework, encompassing the proactive, reactive, and optimization modules. Section \ref{sec: Experiment} compares experimental results from the proposed framework and prevalent frameworks. Section \ref{sec: Deploy} describes the deployment of OptScaler along with the online results at Alipay. Section \ref{sec: Conclude} concludes this paper.

\section{Related Works} \label{sec: Literature}

\begin{table*}[tb]
\caption{Comparison of the proposed OptScaler with the existing representative works on autoscaling}
\begin{center}
\begin{threeparttable}
    \begin{tabular}{llcc}
    \hline
    Type & Literature & Optimization\tnote{c} & Uncertainty\tnote{d} \\
    \hline
    & Zhou2023AHPA\cite{zhou2023ahpa}, Wang2023\cite{Wang2023FullSA}, Das2016\cite{das2016automated} & \ding{55} & \ding{55} \\
    & Dezhabad2018\cite{Dezhabad2018}, Zhang2013\cite{zhang2013dynamic}, Roy2011\cite{roy2011efficient} &  \ding{51} &  \ding{55} \\
    Proactive & Poppe2023\cite{poppe2023proactive}, Qiu2023AWARE\cite{qiu2023aware}, Xue2022\cite{xue2022meta}, Jamshidi2016\cite{jamshidi2016managing} & \ding{55} & \ding{51}  \\
    & Pan2023MagicScaler\cite{pan2023magicscaler}, Qian2022RobustScaler\cite{qian2022robustscaler}, Luo2022\cite{10.1145/3542929.3563477} & \ding{51} & \ding{51}  \\
    \hline
    & Cahoon2022Doppler\cite{cahoon2022doppler}, Rzadca2020Autopilot\cite{49174}, Farokhi2016\cite{farokhi2016hybrid}  & \ding{55} & \ding{55}  \\
    Reactive & Qiu2020FIRM\cite{qiu2020firm} & \ding{55} & \ding{51} \\
    & Gaggero2018\cite{gaggero2018model} & \ding{51}  &  \ding{55} \\
    & Persico2017\cite{persico2017fuzzy} & \ding{51} & \ding{51} \\
    \hline
    Hybrid\tnote{a} & Singh2021RHAS\cite{singh2021rhas}, Jv2018HAS\cite{jv2018has}, Ali2012\cite{ali2012adaptive}, Urgaonkar2008\cite{urgaonkar2008agile} & \ding{55} & \ding{55} \\
    \hline
    Collaborative\tnote{b} & \textbf{OptScaler} & \ding{51} & \ding{51} \\
    \hline
    \end{tabular}
    \begin{tablenotes}
    \footnotesize
    \item[a] The work employs both proactive and reactive modules but they operate independently and may produce conflicting scaling decisions.
    \item[b] The work orchestrates proactive and reactive modules to produce a holistic scaling decision.
    \item[c] The work adopts any optimization technique that can handle practical scaling restrictions and increase interpretability of scaling decisions.
    \item[d] The work considers uncertainty in the cloud system.
    \end{tablenotes}
\end{threeparttable}
\label{tab:literature}
\end{center}
\end{table*}

Autoscaling methods play a critical role in managing cloud resources, with active research in both theoretical and practical implementation across various cloud systems. \autoref{tab:literature} compares our work with representative literature. It is observed that the majority of autoscalers rely on standalone proactive or reactive modules.

In the context of workload prediction, statistical models such as Auto regression or Exponential smoothing are commonly applied \cite{jv2018has, jamshidi2016managing, zhang2013dynamic, qiu2023aware,jin2023timellm}, with recent trends also embracing deep learning models like RNN and Transformer \cite{10.1145/3542929.3563477, pan2023magicscaler,wang2023timemixer,liu2023itransformer}. When it comes to scaling decisions, existing works predominantly use queuing models to make stable decisions. Some also employ analytical methods like PID \cite{9599238, Dynamic_EC2} or pretrained estimators \cite{zhou2023ahpa}, while others rely on customized rules based on prior knowledge \cite{jv2018has, 49174, azure,wang2022end,wang2023flow}. Additionally, optimization techniques and Reinforcement Learning (RL) have also been utilized, and a minority of literature addresses noise through fuzzy decisions \cite{persico2017fuzzy, jamshidi2016managing,wang2024neuralreconciler}.

However, as pointed out by \cite{straesser2022not}, a standalone proactive or reactive autoscaler may fall short of production-ready standards. Only a small number of studies (15 out of 104) have attempted to harness the power of hybrid autoscaling, with the dominant approach being to \textbf{choose between} scaling decisions made by proactive and reactive modules. For instance, \cite{ramperez2021flas, al2017elasticity, ali2012adaptive} adopted the reactive decision when the two conflicted, and \cite{singh2021rhas, jv2018has, urgaonkar2008agile} switched between proactive and reactive decisions based on predefined rules. Notably, our proposed framework revolutionizes this approach by fostering full collaboration between the proactive and reactive modules, with both contributing to the final scaling decision. The former provides future workload inputs, while the latter traces system dynamics between workloads and scaling metrics through real-time feedback.

Previous researchers have also explored leveraging optimization techniques such as Model Predictive Control (MPC) in autoscaling. For example, MPC has been used to adjust the number of virtual machines to meet response time requirements \cite{roy2011efficient}, study the impact of control frequency on resource efficiency and proactive scaling overhead \cite{ogawa2018hierarchical}, solve server placement problems \cite{gaggero2018model, zhang2013dynamic}, and facilitate combined horizontal and vertical scaling with load distribution among machines \cite{incerto2018combined}. However, to our knowledge, MPC (or any other optimization technique) has never been applied in a hybrid/collaborative autoscaler. Our paper advances this area of research by presenting a new approach to exploring the potential of optimization in facilitating collaborative autoscaling.

\section{Preliminaries}
\label{sec: Preliminaries}
In this section, we first present background about cloud resource deployment. Then, we briefly introduce workload patterns and scaling metrics, both are exploited by OptScaler for robust scaling.

\subsection{Basics of Resource Deployment in the Cloud}
In cloud computing, computing resources are organized within a multilevel hierarchical structure. For instance, within Alipay's Recommendation Platform, the cloud computing resource is logically divided into several \textbf{clusters}~\cite{Cluster}. Each cluster encompasses a variable number of \textbf{nodes} (or machines) and supports a group of LRAs. Each node possesses its own dedicated (and usually equalized) computing resources and can be virtually segmented into multiple \textbf{containers}~\cite{Containers}. Each container has the capacity to accommodate a single LRA, and for fault tolerance purposes, an LRA's workload is evenly spread across all host containers. To achieve a well-balanced loading scheme, each node in the cluster deploys one container for each LRA, resulting in all nodes sharing an identical container configuration.

It is important to note that two practical constraints impact the scaling within the LRA scenario. In order to ensure safety, upper and lower bounds for the number of nodes in each cluster must be maintained. Furthermore, the addition or removal of multiple nodes can be conducted in parallel, but there are limitations on the speed and concurrency of node addition/deletion.

\subsection{Patterns of Workload and CPU Utilization} \label{sec: Patterns}
\label{sec: workload_cpu_pattern} 
\begin{figure}[tb]
	\centering
	\includegraphics[width = \columnwidth]{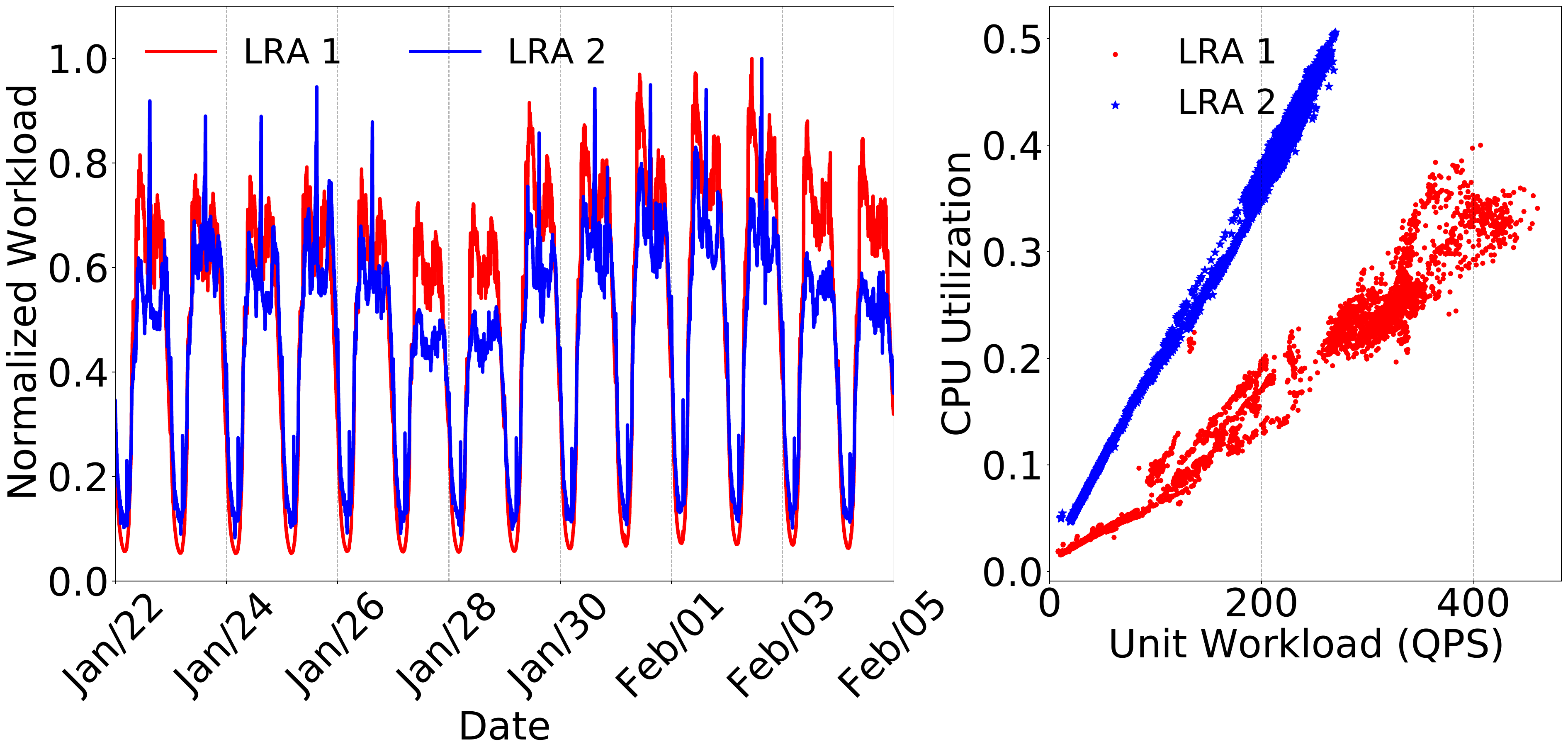}
	\caption{Periodical workloads (left) and linear correlation between unit workload and CPU utilization (right) for two randomly chosen LRAs from Alipay.}
        \Description{}
	\label{fig:patterns_workload_CPU}
\end{figure}

Variations in workload patterns impact the predictive complexity, significantly affecting the performance of the proactive module and consequently, the entire autoscaling framework.
For accurate workload prediction, it is essential that workloads exhibit periodic behaviors. Using randomly sampled workloads of two LRAs from Alipay (from Jan/22/2024 to Feb/04/2024), we have depicted the normalized workload patterns in the left subplot of \autoref{fig:patterns_workload_CPU}. These patterns reveal periodic behaviors in both LRAs, indicating the potential benefits of autoscaling. Specifically, LRA 1 demonstrates daily fluctuations, while LRA 2 exhibits a weekly pattern with lower workloads during weekends.

Moreover, as highlighted in~\cite{SLA_Manage}, there exists a close relationship between SLOs and system metrics such as CPU and memory utilization within a cloud environment, directly influencing the efficacy of cloud resource management. With the primary resource stress stemming from the limited CPU capacity, CPU utilization has been selected as the primary scaling metric. To strike a balance between SLOs and resource conservation, CPU utilization should closely approach (but not exceed) a designated threshold. Consequently, understanding the dynamics between the unit workload (i.e., workload of LRAs in each node) and the resultant CPU utilization is crucial for our framework design.

Initially, we focus on a single node hosting a solitary LRA, assuming that this rule applies similarly when hosting multiple LRAs (i.e., co-location). As shown in the right subplot of \autoref{fig:patterns_workload_CPU}, we observe an explicit linear correlation between the unit workload and the mean CPU utilization. This correlation underscores two key insights: 1) LRAs exhibit varying impacts on CPU utilization. LRA 2 proves to be more resource-intensive than LRA 1, evident from the significantly higher CPU utilization for LRA 2 compared to LRA 1 under the same unit workload (e.g., 200 QPS or queries per second); 2) Under a linear estimator, the variance of residuals sharply increases as the unit workload rises, exacerbating estimator uncertainty. To address these insights, we propose two enhancements to the standard linear estimator. Firstly, we incorporate an uncertainty term, with the variance being a function of the unit workload. Secondly, we implement an online linear regression (OLR) \cite{mohri2018foundations} as a feedback mechanism to dynamically adapt the linear estimator at runtime with the latest feedback, without requiring complete retraining of the estimator with the entire dataset. More details on these enhancements will be provided in Section \ref{sec: cpu_estimator}.

\section{Proposed Framework} \label{sec: Proposed}
\begin{figure}[tb]
	\centering
	\includegraphics[width = 0.8\columnwidth]{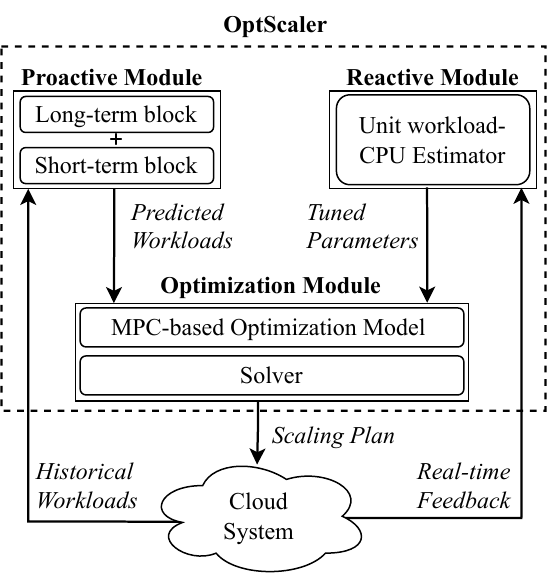}
	\caption{Flowchart of the proposed OptScaler framework.}
        \Description{}
	\label{fig:flowchart}
\end{figure}

When implementing horizontal scaling, OptScaler could focus on the scaling of the number of nodes in each cluster for making holistic scaling decisions. The OptScaler framework, illustrated in \autoref{fig:flowchart}, comprises the following components:

\begin{itemize}
    \item \textbf {Proactive module} consists of a workload prediction model trained on historical workloads, providing multi-timestep predictions of future workloads for all LRAs of interest;
    \item \textbf {Reactive module} offers a self-tuning estimator of CPU utilization. It is initialized using historical data of the unit workloads and the corresponding CPU utilization, and continuously fine-tuned by real-time system feedback;
    \item \textbf{Optimization module} takes input from both proactive and reactive modules, producing an optimal scaling plan under practical restrictions. The plan is then deployed to the cloud system.
\end{itemize}

\subsection{Proactive Module}
\label{sec: ts_model}

In Proactive Module, we employ a workload prediction model to anticipate the future workloads of different LRAs.
The casting results are further forwarded to the optimization module for CPU utilization management.
The following challenges must be addressed for effective autoscaling: 
1) Long-term forecasting: In order to complete node addition/removal within each time interval amidst rapid workload fluctuations, the optimization module must plan several steps in advance. This necessitates long-term forecasting, such as predicting workload a day in advance.
2) Model efficiency: Due to the presence of multiple LRAs with diverse weekly and daily workload patterns, constructing a separate prediction model for each LRA would be inefficient. It is advantageous to have a single model that can accommodate the various workload patterns, providing convenience for the online service management.

\begin{figure}[tb]
    \centering
    \includegraphics[width=\columnwidth]{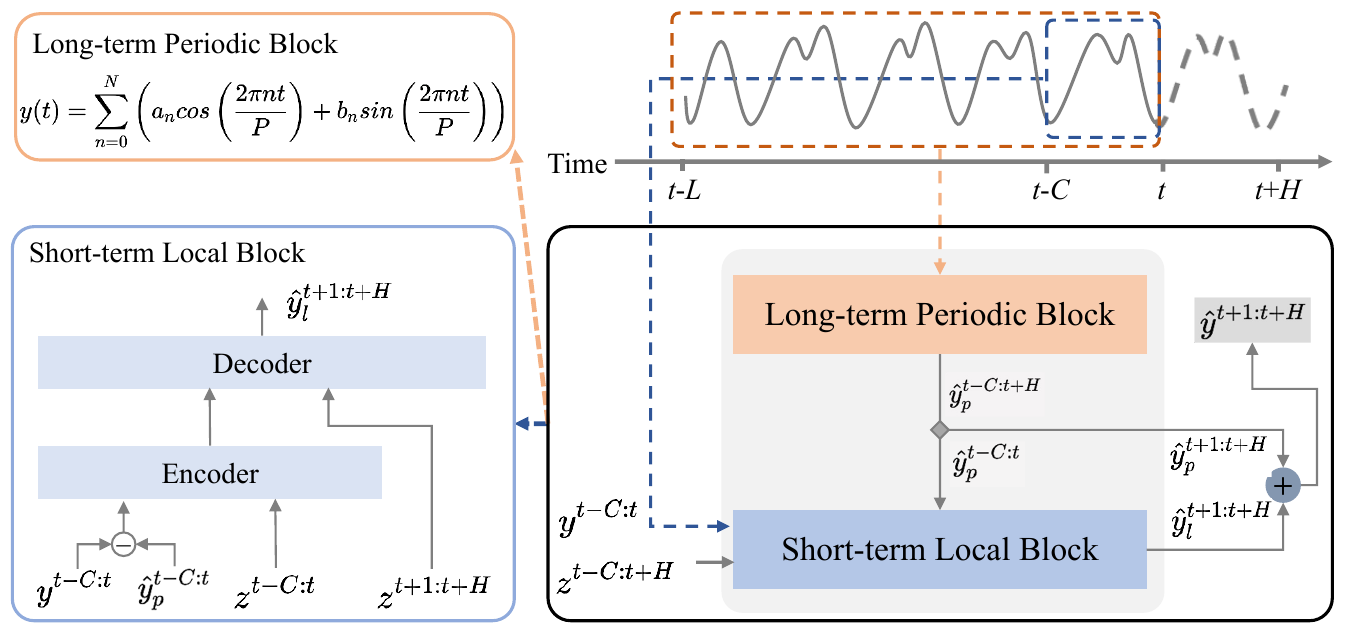}
    \caption{Framework of the workload prediction model in the Proactive Module.}
    \Description{}
    \label{fig:ts_model}
\end{figure}

To address the above challenges, a workload prediction model is introduced in \autoref{fig:ts_model}.
Formally, let $y_n^t \in \mathbb{R}$ denotes the workload at time step $t$ of the $n$-th LRA, the task is to predict the future values $\boldsymbol{y}_n^{t+1:t+H}=[y_n^{t+1},\ldots,y_n^{t+H}]$ based on the historical values $\boldsymbol{y}_n^{t-C:t}=[y_n^{t-C},\ldots,y_n^{t}]$ and other covariates. 
We herein use bold symbols to denote \textbf{vectors}.
Typically, it can be formalized as:
\begin{equation}
\label{eq:forecasting}
\boldsymbol{y}_n^{t+1:t+H} = \mathcal{F}_{\Theta}(\boldsymbol{y}_n^{t-C:t},\boldsymbol{z}_n^{t-C:t+H},n),
\end{equation}
where $\mathcal{F}_\Theta(\cdot)$ denotes the workload prediction model, $\Theta$ is the learning parameters, $\boldsymbol{z}_n^{t-C:t+H}$ is the known covariates (e.g., date) of the $n$-th LRA, and $n \in \{1,2,\ldots,N\}$ is the index of LRA\footnote{For simplicity, we ignore the index $n$ of LRA in the following.}.

It is evident that $\mathcal{F}_\Theta(\cdot)$ handles the workload prediction task for various LRAs with the parameter $\Theta$ and distinguishes LRAs with index $n$. Specifically, $\mathcal{F}_\Theta(\cdot)$ comprises a Long-term Periodic Block and a Short-term Local Block to accomplish this task effectively.
\begin{itemize}
    \item \textbf{Long-term Periodic Block}. 
    Capture too long series input for forecasting workloads of LRAs, each with varying weekly or daily patterns, can be challenging.
    To address this, we employ the Fourier series to represent the diverse periodicities of different LRAs. 
    With a periodicity $P$ and truncation order $\hat{N}$ of the Fourier series expansion, 
    the long-term periodicity can be formalized as $y^t =\sum_{n=0}^{\hat{N}}(a_ncos(\frac{2\pi nt}{P})+b_n sin(\frac{2 \pi nt}{P}))$.
    The Fourier coefficients $[a_0,b_0,\ldots,a_{\hat{N}},b_{\hat{N}}]$, learned from historical data, facilitate the inference of future periodicity without the need for extensive historical inputs.
    
    \item \textbf{Short-term Local Block}.
    This block learns the local pattern of workloads, such as local trend and influence.
    For this purpose, we apply a Seq2seq structure consisting of three steps:
    1) According to historical values $\boldsymbol{y}^{t-C:t}$ and their periodic estimation $\boldsymbol{\hat y}_p^{t-C:t}$, we calculate the residual of historical series $\boldsymbol{y}_l^{t-C:t}=\boldsymbol{y}^{t-C:t}-\boldsymbol{\hat y}_p^{t-C:t}$;
    2) We extract local pattern of historical residual with series $\boldsymbol{y}_l^{t-C:t}$ and corresponding covariates $\boldsymbol{z}^{t-C:t}$. A linear-complexity Flow-Attention~\cite{wu2022flowformer} $f(\boldsymbol{q},\boldsymbol{k},\boldsymbol{v};\boldsymbol{\theta}_{attn})$ is applied for long sequences with high efficiency. Specifically, let $\boldsymbol{h}_{in}=[\boldsymbol{y}_l^{t-C:t}, \boldsymbol{z}^{t-C:t}]$ denote the whole historical representation, 
    the encoder is $\boldsymbol{h}_{enc}=f(\boldsymbol{h}_{in},\boldsymbol{h}_{in},\boldsymbol{h}_{in};\boldsymbol{\theta}_{enc})$, where $\boldsymbol{\theta}_{enc}$ is the parameters of encoder;
    3) Based on the encoder output $\boldsymbol{h}_{enc}$ and covariates $\boldsymbol{z}^{t-C:t+H}$, the decoder estimates the local residual in future. With the same Flow-Attention structure, the decoder is $\boldsymbol{\hat y}_l^{t+1:t+H}=f(\boldsymbol{z}^{t-C:t}, \boldsymbol{z}^{t+1:t+H}, \boldsymbol{h}_{enc}; \boldsymbol{\theta}_{dec})$, where $\boldsymbol{\theta}_{dec}$ is the learning parameter.
\end{itemize}

In the Short-term Local Block, two practical techniques are employed to enhance long-term forecasting. 
The first technique, Flow-Attention, demonstrates computational efficiency with linear complexity  relative to the length of the series.
The second technique involves the use of a non-autoregressive decoder, where the forecasting of time step $t$ is independent of the results from the preceding time step. This non-autoregressive decoder helps avoid accumulative errors and offers swift inference speed.

As illustrated in the flowchart Figure~\ref{fig:flowchart}, the predicted workload of each LRA is the sum of the outputs from Long-term Periodic Block and Short-term Local Block:
\begin{equation}
\boldsymbol{\hat y}^{t+1:t+H}=\boldsymbol{\hat y}_p^{t+1:t+H}+\boldsymbol{\hat y}_l^{t+1:t+H}
\end{equation}
During model training, quantile loss \cite{park2022learning,zhu2021mixseq} helps minimize the ratio of the forecasting results that are lower than the actual values. This approach enhances the reliability of subsequent scaling decisions from the perspective of proactive module.

\subsection{Reactive Module}
\label{sec: cpu_estimator}

\begin{algorithm}[tbp]
\caption{Widrow-Hoff algorithm in the Reactive Module}
\label{alg:olr}
\begin{algorithmic}
    \State \textbf{Input:} Choose parameter $\eta > 0$
    \State \textbf{Initialize:} $\boldsymbol{w_k}^0$ according to historical data training
    \For{$t=1$ to $T$}
        \State Get unit workload $\frac{\boldsymbol{y}^t}{x^t}$ $ \in \mathbb{R}^N$
        \State Estimate $\hat{c}^t = f(\frac{\boldsymbol{y^t}}{x^t}) \in \mathbb{R}$
        \State Observe $c^t \in \mathbb{R}$
        \State Get feedback error $e = \hat{c}^t - c^t$
        \State Update $\boldsymbol{w_k}^{t} = \boldsymbol{w_k}^{t-1} - \eta e \frac{\boldsymbol{y^t}}{x^t}$
    \EndFor \\
    \Return $\boldsymbol{w_k}^{T}$
\end{algorithmic}
\end{algorithm}

In Reactive Module, we build a linear estimator with an uncertainty term to map the unit workloads of LRAs to the average CPU utilization of all nodes in the cluster, which could be formulated as:
\begin{equation}
c = f\left(\frac{\boldsymbol{y}}{x} \right) = w_b + \frac{\boldsymbol{y}^\top \boldsymbol{w_k}}{x}  + \epsilon 
\label{eq:cup_util_1}
\end{equation}
where $f(\cdot)$ denotes the estimator;
$\boldsymbol{y} \in \mathbb{R}^N $ is the workload vector for all $N$ LRAs;
$x$ is the number of nodes in the cluster to be decided;
$c$ denotes the estimated CPU utilization, which is a function of unit workload $\frac{\boldsymbol{y}}{x}$; 
$\boldsymbol{w_k}$ and $w_b$ are the slope and intercept for the linear model, respectively. 
$\boldsymbol{w_k}$ can be seen as the weights of the LRAs, and $w_b$ is naturally the overhead load of CPU; $\epsilon$ is the uncertainty term compensating for the residual of a linear model. $\epsilon$ obeys the following normal distribution:
\begin{equation}
\epsilon \sim N (0, \sigma^2(\frac{\boldsymbol y}{x}) ) \sim N (0, (\sigma_b + \frac{\boldsymbol y^\top \boldsymbol{\sigma}_k}{x})^2 ) 
\label{eq:cup_util_2}
\end{equation}
where the standard deviation of $\epsilon$ is modelled as a linear function (with coefficient $\boldsymbol{\sigma}_k$ and $\sigma_b$) of unit workload. In this way, we ensure that the uncertainty term obeys the rules shown in \autoref{fig:patterns_workload_CPU} that the variance of CPU utilization increases with larger unit workloads.

Based on the above formulations, we now introduce the reactive mechanism. 
Initially, all the parameters (including $\boldsymbol{w_k}$, $w_b$, $\boldsymbol{\sigma}_k$ and $\sigma_b$) are initialized using maximum likelihood estimation ~\cite{myung2003tutorial} on the historical data with non-negative constraints. 
Subsequently, we employ OLR~\cite{mohri2018foundations} as a means to adjust $\boldsymbol{w_k}$ using real-time feedback.
The feedback mechanism of OLR resembles that of a supervised learning algorithm, aiming to minimize the cumulative square loss of a linear function in an online setting. 
The pseudocode is shown in Algorithm \ref{alg:olr}, where $\eta$ is a parameter to control the strength of feedback tuning, in analogy with the learning rate in a machine learning context, and the real-time feedback error $e = \hat{c}^t - c^t$ is derived from the cloud.
In practice, each time when new feedback is available, we could choose the most recent (i.e., $T = 1$) or a series of (i.e., $T > 1$) feedback, and apply OLR to update $\boldsymbol{w_k}$ quickly. Instead of retraining the linear model using the whole dataset each time, OLR focuses on the newest $T$ feedback, and its simplicity and efficiency become key advantages in an online system.

\subsection{Optimization Module}
Our optimization module is designed to enable collaboration between the proactive module and the reactive module, fostering holistic decision-making. It leverages MPC to dynamically take inputs from both modules and consider all the practical constraints in the cloud. As stated in Section \ref{sec: Intro}, MPC shows promise for addressing future workload bursts when the scaling speed is restricted.

In general, for each control action, MPC takes the latest system state and solves a constrained optimization model over a sliding window of future time intervals. It only applies the first solution over the horizon and repeats the procedure the next time \cite{holkar2010overview}. 
In our context, we repeatedly optimize scaling decisions of future $D$ time intervals at the beginning of every time interval.
Each time interval $d \in \{1,2,\ldots, D\}$ spans $h$ minutes (e.g., 30 minutes), which can be much longer than the resolution of workload predictions (e.g., 1 minute). Given that each cluster is highly autonomous with its exclusive LRAs and resources, we will focus on a single cluster to build the optimization model. At time step $t$ (i.e., the beginning of a time interval), the following model of $D$ time intervals is solved:
\begin{align}
\max_{\boldsymbol{u}} \quad & \sum_{d=1}^D c^d  \label{obj}\\ 
s.t. \quad & |u^d|\leq \frac{h}{\tau} \cdot s, &\forall d  \label{cons:speed}\\
& x^d = x^0 + \sum_{j=1}^{d}u^j, &\forall d \label{cons:state}\\
& c^d= \max_{j \in \{d,d+1\}} f\left(\frac{\boldsymbol{y}^j}{x^d}\right), &\forall d \label{cons:maxcpu}\\
& c^d\leq c^*, &\forall d \label{cons:cpubound}\\
& X^{min} \le x^d\le X^{max}. &\forall d \label{cons:xbound}
\end{align}

\noindent
\textbf{Inputs}: 
$x^0$ represents the initial number of nodes, capturing the number of nodes deployed at time $t$. 
$\boldsymbol{y}^d$ denotes the peak value of the predicted workload in time interval $d$ (i.e., $\boldsymbol{y}^{t+(d-1)h:t+dh}$) taken from the proactive module. $f(\cdot)$ denotes the CPU utilization estimator taken from the reactive module.

\noindent
\textbf{Variables}: $u^d$ denotes the maximum change of the number of nodes in time interval $d$, and $\boldsymbol{u}=[u^1,\ldots,u^D]\in \mathbb{R}^{D}$ is the vector. Only $u^1$ is returned for deployment. $x^d$ denotes the number of nodes provisioned during time interval $d$.

\noindent
\textbf{Objective}: \eqref{obj} is the objective that maximizes the CPU utilization $c^d$ (without exceeding a desired level $c^*$) over all $D$ time intervals.

\noindent
\textbf{Constraints}: \eqref{cons:speed}-\eqref{cons:xbound} are explained as follows: 
\begin{itemize}
    \item Constraint \eqref{cons:speed}: The speed limitation of scaling, where $h$ is the time interval between two successive scaling actions, while $\tau$ and $s$ represent the time cost and the concurrency of single node addition/removal action, respectively;
    \item Constraint \eqref{cons:state}: The state transition equation that defines the relationship between $x^d$ and $x^0$;
    \item Constraint \eqref{cons:maxcpu}: 
    The calculation of CPU utilization $c^d$, estimated as the maximum value in two adjacent time intervals \{$d, d+1$\}. This is due to the fact that scaling decisions are made in advance, where the resource for time interval $d+1$ is provisioned ahead during $d$, leading to $f(\boldsymbol{y}^{d+1}/x^d)$); 
    \item Constraint \eqref{cons:cpubound}: The upper bound $c^*$ for $c^d$ ;
    \item Constraint \eqref{cons:xbound}: The upper and lower bounds for the number of nodes $x^d$.
\end{itemize} 

\noindent
\textbf{Reformulation}: Due to uncertainty term $\epsilon$ in $f(\cdot)$ as shown in \eqref{eq:cup_util_1}, constraint \eqref{cons:cpubound} is easily violated. We reformulate constraints \eqref{cons:maxcpu} and \eqref{cons:cpubound} by chance constraint \eqref{cons:cc} that guarantees satisfaction of constraint \eqref{cons:cpubound} with probability $\alpha$: 
\begin{equation} \label{cons:cc}
 \mathbb{P} \{ c^* \ge c^d =\max_{j\in \{d,d+1\}} f\left(\frac{\boldsymbol{y}^j}{x^d}\right) \} \ge \alpha,  \quad \forall d
\end{equation} 

\noindent
\textbf{Equivalent transformation}: Consider the formulation of $f(\cdot)$ in \eqref{eq:cup_util_1}, constraint \eqref{cons:cc} can be converted to a deterministic equivalent form under the theory of chance-constraint method \cite{shapiro2021lectures}:
\begin{align}
\label{cons:cc_deter}
c^* \ge c^d & = \frac{1}{x^d} \max_{j\in \{d,d+1\}} (\psi^{-1}(\alpha) \cdot {\boldsymbol {\sigma}_{k}} + {\boldsymbol{w}_{k}})^\top \boldsymbol{y}^j \nonumber \\
 & + w_{b} + \psi^{-1}(\alpha) \cdot\sigma_{b}, \hspace{2.22cm} \forall d 
\end{align}
where $\psi^{-1}(\cdot)$ is the inverse of the cumulative distribution function of a standard normal distribution. For example, $\psi^{-1}(\alpha) \approx 1.28$ when $\alpha=0.9$. Define $m^d = \max_{j\in \{d,d+1\}} (\psi^{-1}(\alpha)\cdot {\boldsymbol {\sigma}_{k}} +  {\boldsymbol{w}_{k}})^\top \boldsymbol {y}^j$, which can be calculated before solving the optimization model, \eqref{cons:cc_deter} can be easily rearranged to an equivalent linear form:
\begin{align}
\label{cons:cc_linear}
x^d \ge \frac{m^d}{c^* - w_{b}  - \psi^{-1}(\alpha) \cdot\sigma_{b}}. \hspace{1.4cm} \forall d 
\end{align}
Constraint \eqref{cons:cc_deter} also restricts $c^d$ from exceeding $c^*$, hence, by omitting constant terms, objective \eqref{obj} is equivalent to
\begin{equation} \label{obj:new}
 \max_{\boldsymbol{u}} \quad  \sum_{d=1}^D \frac{m^d}{x^d}.
\end{equation}

\noindent
\textbf{Solver}: Now the optimization problem \eqref{obj}-\eqref{cons:xbound} under uncertainty can be reformulated as a robust counterpart with a non-linear objective \eqref{obj:new} and linear constraints \eqref{cons:speed}, \eqref{cons:state}, \eqref{cons:xbound}, \eqref{cons:cc_linear}.  This reformulated problem can be solved using off-the-shelf solvers (e.g. IPOPT \cite{wachter2006}) and the solution $u^1$ can be rounded for implementation.

\noindent
\textbf{Robustness}:  In summary, we utilize MPC and chance constraints to ensure the satisfaction of practical constraints and bolster the robustness of the solution.

\section{Numerical Results} \label{sec: Experiment}
In this section, we provide a detailed exploration of the three research questions in Section \ref{sec: Intro}. We specifically investigate the performance of the proposed OptScaler in the context of co-located LRAs with diverse workload patterns. Question 1 is addressed in Section \ref{exp:fcst_model}, while Questions 2 and 3 are discussed in Section \ref{sec: compare_framework}.

\subsection{Workload Dataset} \label{sec: data}
To facilitate a comparison of different autoscalers, we conduct experiments using datasets from the public Azure Functions Trace~\cite{254430} and the internal Alipay Applets Trace. We select five workloads from LRAs in the Azure Trace and three from the Alipay Trace. Each dataset spans a duration of two weeks (Jul/15 to Jul/28) and maintains a minute-level resolution. Each LRA is identified by the last four bits of its hash ID.
We categorize the workloads into three sets, denoted as $S_A$, $S_B$, and $S_C$, based on the complexity of their workload patterns. Sets $S_A$ and $S_C$ are derived from the Azure Trace, while set $S_B$ is derived from the Alipay Trace. LRAs within each set will be co-located in a cluster, resulting in three clusters of interest. Subsequently, experiments on workload prediction and autoscaling will be carried out for these three clusters.

\begin{figure}[tb]
	\centering
	\includegraphics[width = \columnwidth]{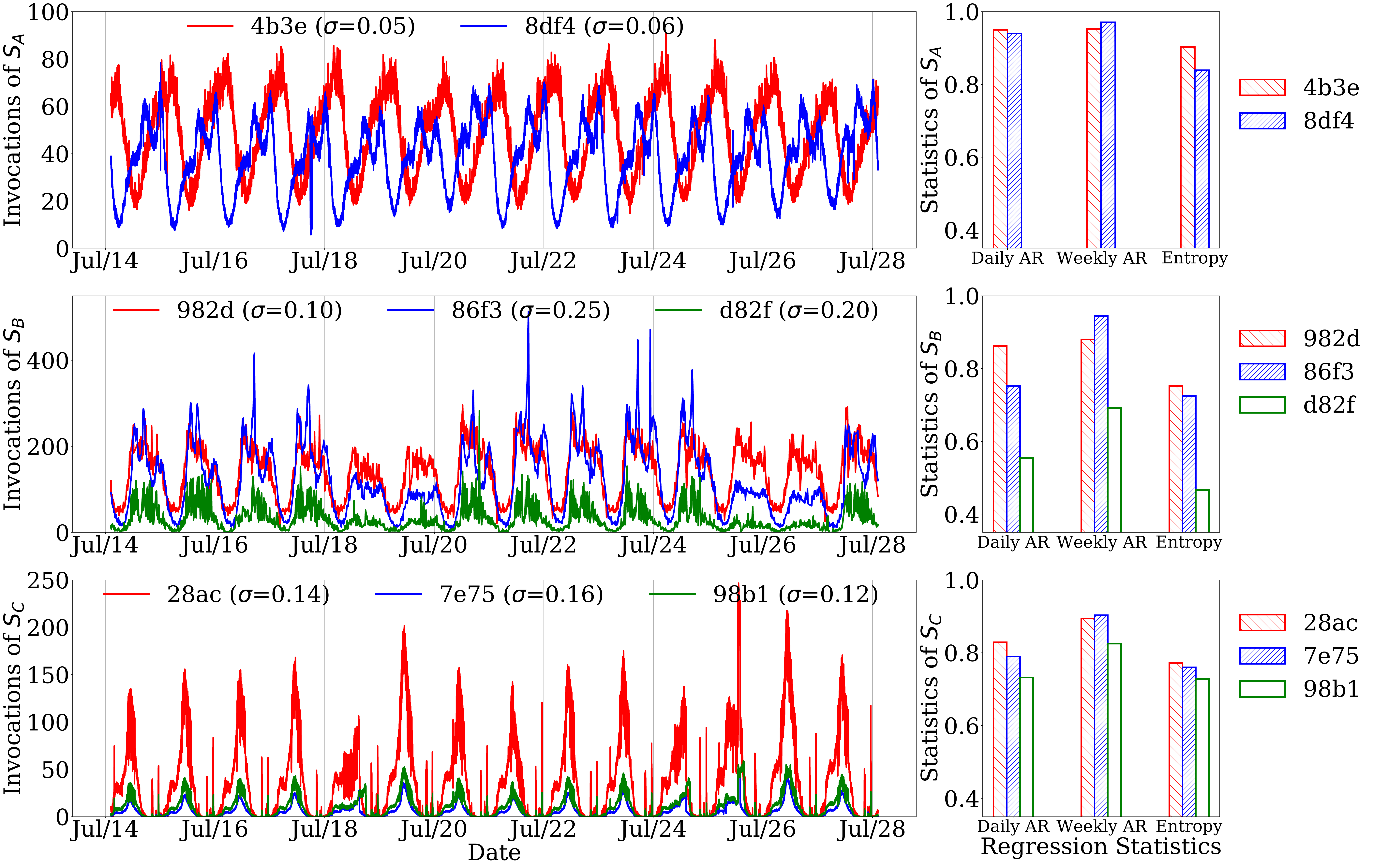}
	\caption{Time series (left) and regression statistics (right) for Azure Trace ($S_A$, $S_C$) and Alipay Trace ($S_B$). The complexity of prediction and scaling increases from $S_A$ to $S_B$ and $S_C$.}
        \Description{}
	\label{fig:Azure_series_new}
\end{figure}

\autoref{fig:Azure_series_new} provides details about workloads in the three sets. The left three subplots show the time series of workloads (measured in QPS). Along with the legend, $\sigma$ denotes the standard deviation of daily peaks of the normalized series.
A higher $\sigma$ value signifies a more pronounced variation in the daily peaks, posing a greater challenge for prediction models. The right three subplots show the regression statistics: daily coefficient of auto-correlation (i.e., \textbf{Daily AR}), weekly coefficient of auto-correlation (i.e., \textbf{Weekly AR}), and spectral entropy (i.e., \textbf{Entropy}); the first two statistics detect the strength of daily/weekly periodicity, and the spectral entropy normalizes the overall complexity of the workloads \cite{Jiang_KATS_2022}. We take the bar value of Entropy as 1.0 minus the original spectral entropy, such that for all the three regression statistics, a lower bar indicates a higher level of prediction difficulty. Clearly, the complexity for $S_A$ is much lower than that of $S_B$ and $S_C$.

In set $S_A$, both workloads (4b3e and 8df4) exhibit steady fluctuations of daily peaks with a low standard deviation ($\sigma \le 0.1$) and demonstrate clear daily periodicity, with both daily and weekly auto-correlation values exceeding $0.9$.
For set $S_B$, the workloads display weak periodicity (e.g., d82f), elusive spikes (e.g., 86f3), and volatile daily peaks, with the standard deviation reaching up to $0.25$.
In set $S_C$, three workloads (28ac, 7e75, and 98b1) show similar trends. Specifically, they all experience dramatic workload increases and unstable daily peaks, with a standard deviation of $\ge 0.1$. Additionally, two spikes of different magnitudes between 6 PM and 9 PM each day further complicates accurate prediction and robust scaling, especially for workload 28ac.
In summary, the workload patterns in the three sets effectively represent various scenarios of co-located LRAs, each with distinct levels of complexity. These diverse patterns present a significant challenge and effectively distinguish our proposed autoscaling framework from others.

\subsection{Experiment Settings}
During offline simulations, we emulate the actual implementation of OptScaler in a production cloud. The data from the initial 12 days will be used as training data for the workload prediction model. In the last 2 days (Jul/27 to Jul/28), OptScaler will be triggered every half hour (i.e., $h=30$) to control the CPU utilization to approach the target $c^* = 0.5$.

Additionally, we take into account the scaling limitations, determined by the time interval $h$, the time cost $\tau$ of a single node addition/removal, and the concurrency of scaling $s$ in \eqref{cons:speed}. We set $\tau=5$ minutes and $s=4$ as in a real system. Consequently, with $h=30$, a maximum of 24 nodes could be added or removed in each time interval.
Furthermore, we establish the upper bound $X^{max} = 400$, and the lower bound $X^{min} = 80$ for $S_A$ and $S_B$, and $X^{min} = 20$ for $S_C$, considering that the total workload often approaches zero in the latter case. We set the parameter $\alpha=0.95$ for the chance constraint \eqref{cons:cc} for all optimization-based frameworks. As for OLR (see Algorithm \ref{alg:olr}), we select the feedback strength $\eta = 2e-4$ by minimizing the cumulative loss on the training data.

\subsection{Comparison of Workload Prediction Models}
\label{exp:fcst_model}
\subsubsection{Prediction Models}
To address Question 1, we compare our workload prediction model with several traditional and state-of-the-art time series prediction models, namely ARIMA~\cite{makridakis1997arma}, Autoformer~\cite{wu2021autoformer}, DEPTS~\cite{fan2022depts}, and NBEATS~\cite{oreshkin2019n}.
ARIMA represents a traditional time series prediction model in statistics. Due to the workload pattern in Section~\ref{sec: workload_cpu_pattern}, we adopt ARIMA with the exogenous regressors of Fourier series.
Autoformer is a neural network model designed for long-term prediction based on auto-correlation mechanism and a seq2seq structure.
DEPTS and NBEATS are two types of neural networks tailored for uni-variate time series prediction, incorporating the neural expansion analysis method.
In the rolling window experiment, we train the ARIMA model every half hour using the most recent 24-hour data, owing to its fast training speed. As for the neural network-based models, training occurs at midnight each day. Throughout the training process, data preceding the training time is included in the training dataset, while the data from the last two days is utilized for validation. All these neural network-based models are trained with the $0.5$-quantile loss.

\subsubsection{Evaluation Metrics}
We focus on the long-term prediction performance, that is, predicting the workload in the following $360$ minutes every $h=30$ minutes.
The metric, weighted absolute percentage error (\textbf{WAPE}), is adopted to evaluate different models:
\begin{align}
    WAPE & = \frac{1}{N} \sum_{n=1}^N \frac{\sum_{(t,\delta) \in \Omega} |y_n^{t+\delta} - \hat{y}_n^{t+\delta}|}{\sum_{(t,\delta) \in \Omega} |y_n^{t+\delta}|},
\end{align}
where $n$ is the index of time series (LRAs in our experiments).
$t$ is the starting time of each rolling window,
and $\delta$ is the increment in time, 
that is, $(t,\delta)$ means the $\delta$-th step after the starting time $t$ of the rolling window.
$\Omega$ denotes the whole evaluation space.

We evaluate the prediction accuracy in two aspects:
1) The \textbf{Minute-level} accuracy, that is, whether the model can accurately predict the workload of each LRA in each minute, since all the models are built on the data in minute-level resolution;
2) The \textbf{Interval-peak} accuracy, that is, whether the model can correctly capture the peak workload in the given interval used in MPC-based optimization, i.e., $\boldsymbol{y}^j \in \mathbb{R}^N$ in constraint~\eqref{cons:maxcpu}.
This accuracy is critical to the subsequent optimization.

\subsubsection{Experimental Results for Question 1} \label{sec: Q1}

\begin{table}[tb]
    \caption{
    Comparison of five prediction models, testing on WAPE for the Minute-level and Interval-peak accuracy. Lower values are better, and the best are bolded.}
    \begin{center}
    \begin{tabular}{lcccccc}
        \toprule
        & \multicolumn{3}{c}{\textbf{Minute-level}} & \multicolumn{3}{c}{\textbf{Interval-peak}}\\
        \cmidrule(r){2-4} \cmidrule(r){5-7}
        & $S_A$ & $S_B$ & $S_C$ & $S_A$ & $S_B$ & $S_C$ \\
        \midrule
        ARIMA          & 0.087 & 0.280 & 0.432 & 0.107 & 0.286 & 0.420 \\
        Autoformer     & 0.079 & 0.258 & 0.467 & 0.085 & 0.274 & 0.454 \\
        DEPTS          & 0.066 & 0.243 & 0.356 & 0.082 & 0.257 & 0.342 \\
        NBEATS         & 0.067 & \textbf{0.218} & 0.352 & 0.085 & 0.178 & 0.377 \\
        Our model      & \textbf{0.063} & 0.224 & \textbf{0.309} & \textbf{0.080} & \textbf{0.171} & \textbf{0.341} \\
        \bottomrule
    \end{tabular}
    \label{tab:ts_pred_res}
    \end{center}
\end{table}

\autoref{tab:ts_pred_res} illustrates the superior performance of OptScaler compared to four other prediction models across LRAs in $S_A$, $S_B$, and $S_C$, as measured by the Weighted Absolute Percentage Error (WAPE) from two different perspectives.
Notably, all methods exhibit considerably better performance on set $S_A$ in comparison to $S_B$ and $S_C$. This observation aligns with the discussion in \autoref{sec: data}, where it is evident that workloads in $S_A$ display relatively strong periodicity and maintain a steady peak trend, whereas those in $S_B$ and $S_C$ demonstrate diverse fluctuations.
Furthermore, our model consistently outperforms ARIMA by a significant margin. This further attests to the effectiveness of our Short-term Local Block, as both models extract the periodic property using Fourier series.
In terms of Minute-level accuracy, our model outperforms the others in $S_A$ and $S_C$ and demonstrates the second best performance in $S_B$. Additionally, in the context of Interval-peak accuracy, a critical aspect in our proposed autoscaling framework, our model outshines all compared models.

\begin{figure}[tb]
    \centering
    \includegraphics[width = \columnwidth]{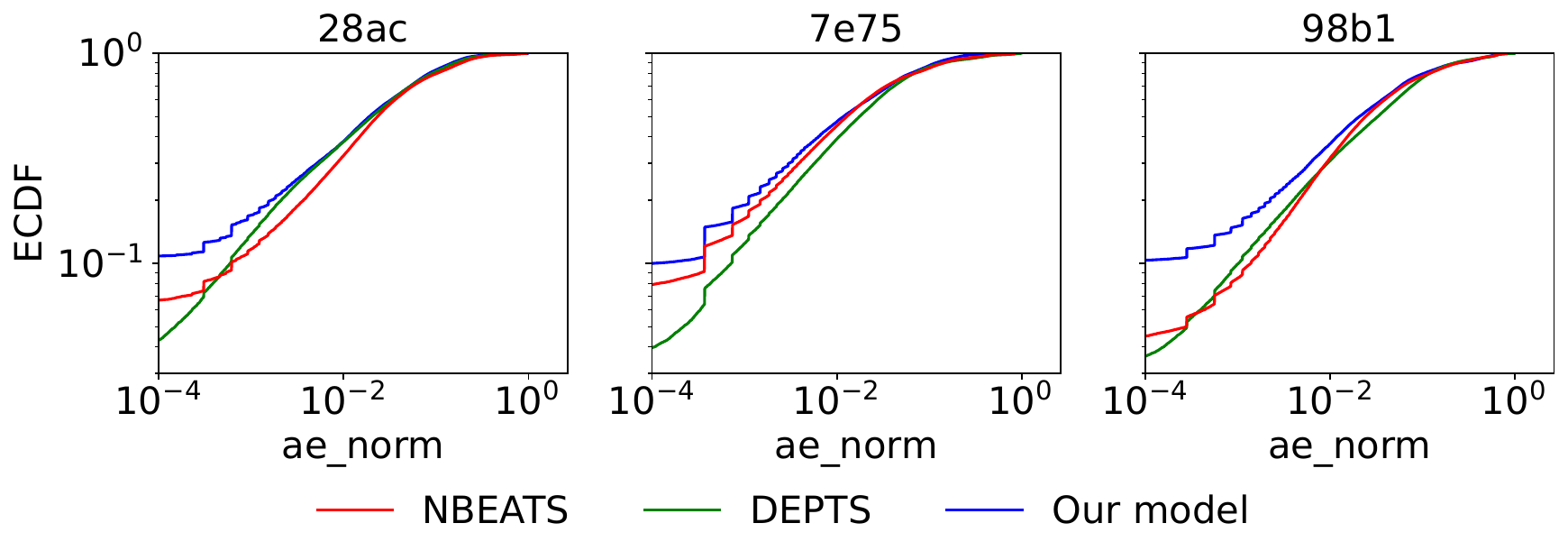}
    \caption{The empirical cumulative distribution function (ECDF) of normalized absolute error ($ae\_norm$) with respect to Minute-level accuracy for three LRAs in set $S_C$. Both x-axis and y-axis are in log scale.}
    \Description{}
    \label{fig:ecdf_sb}
\end{figure}

To further highlight this superiority, we present \autoref{fig:ecdf_sb}, which illustrates the performance of our model in comparison to the two most competitive neural network-based models (NBEATS and DEPTS) presented in \autoref{tab:ts_pred_res}.
Specifically, the Minute-level accuracy of three LRAs in set $S_C$ is compared using their empirical cumulative distribution functions (ECDF) of normalized absolute error ($ae\_norm$, i.e., the absolute error divided by the maximal absolute error for each workload). As depicted in \autoref{fig:ecdf_sb}, our method achieves notably higher empirical cumulative distribution on smaller $ae\_norm$ values. The ECDF curves of our method consistently surpass the others across all LRAs, indicating the superior performance of our method on complex workloads in set $S_C$.

In general, in our model, the Fourier series in long-term periodic block could capture various workload patterns without the necessity of long series input to the model, and the linear-complexity attention mechanism in short-term local block learns the local trend and fluctuation; the combination of these two blocks gives our model an edge over neural network based models like NBEATS and DEPTS in terms of efficiency and accuracy. Our model can also handle multiple time series by one model with fewer parameters, which is much more suitable for the scenario of co-located LRAs.


\subsection{Comparison of Autoscaling Frameworks} \label{sec: compare_framework}
To answer Questions 2 and 3 stated at the beginning of Section \ref{sec: Experiment}, we conduct an ablation experiment to compare the performance of OptScaler with prevalent autoscaling frameworks. Specifically, we respond to Question 2 in Section \ref{sec: Q2} by comparing all four frameworks, given the same input of the prediction results; then, we address Question 3 in Section \ref{sec: Q3} by comparing OptScaler with the best combination from the existing prediction models and autoscaling frameworks.
\subsubsection{Autoscaling Frameworks} \label{sec: Autoscaling_Frameworks}
We investigate four autoscaling frameworks, i.e., reactive, proactive, hybrid, and collaborative:
\begin{itemize}
    \item \textbf{Autopilot} \cite{49174} serves as the representative of \textbf{reactive} autoscalers. Introduced by Google in 2020, it is considered as an industrial benchmark in this context. Autopilot determines the optimal resource configuration by identifying the best-matched historical time intervals to the current window. We implement Autopilot's horizontal scaling based on its publicly available paper \cite{49174};
    \item \textbf{Madu} \cite{10.1145/3542929.3563477} exemplifies a typical \textbf{proactive} autoscaler, leveraging the capacity of both prediction and optimization. We adapt Madu to our context, formulating the optimization model as \eqref{obj}-\eqref{cons:xbound}. As a pure proactive autoscaler, the CPU utilization estimator will not be updated after initialization;
    \item \textbf{HAS} \cite{jv2018has} represents a \textbf{hybrid} autoscaler that adopts the mainstream concept of independent proactive and reactive modules. Without compromising on quality, we integrate HAS into our context, where the reactive decision replaces the proactive one when the observed CPU utilization exceeds a predefined bound. The current bound is set to be $0.9c^*$. The reactive decision at current time $t$ is calculated by $u^t=\max\{X^{min}-x^{t-1},\ \min\{(c^{t-1}/c^* -1) x^{t-1}, \ h\cdot s/\tau, \ X^{max} - x^{t-1}\}\}$, adjusting the resource to eliminate the scaling error at $t-1$. The $\min$ operation indicates $c^t \le c^*$ if workload at $t$ does not increase, while the $\max$ operation ensures the lower bound in \eqref{cons:xbound} holds.
    \item \textbf{OptScaler} is our proposed \textbf{collaborative} framework.
\end{itemize}

Accordingly, we customize the following two parameters for the above autoscaling frameworks: 
\begin{itemize}
    \item $S$ is a parameter for Autopilot that denotes the percentile of the most recent samples of CPU utilization. A higher value of $S$ indicates a more conservative scaling preference. Here we consider two different values for $S$, namely $S \in \{90\%, 95\%\}$ for Autopilot;
    
    \item $D$ represents time horizon in the optimization model for Madu, HAS, and OptScaler. In our context, we use a $30-$ minute interval (i.e., $h=30$) for each $d$. A value of $D = 1$ signifies a standard optimization model focusing on immediate scaling needs, while $D = 11$ represents an MPC-based optimization model that can anticipate future events within the following 6-hour (as per constraint~\eqref{cons:maxcpu}, we need to consider one more time interval than $D$). Therefore, in the case of $D = 11$, we will end our experiment at 6 PM on Jul/28.
\end{itemize}

\subsubsection{Evaluation Metrics}
\label{sec: Metrics}
As an end-to-end autoscaling framework, our primary objective is to achieve a balance between the cost of computing resources and the stringent SLOs in cloud service. With this goal in mind, we have developed three evaluation metrics, where lower metric values indicate better performance:

\begin{itemize}
    \item \textbf{$S_{vr}$} denotes the SLO violation rate, representing the percentage of minutes with SLO violation (i.e., $c^d > c^*$) during the experiment period. Therefore, a high value of $S_{vr}$ indicates a prolonged duration of SLO failure;
    \item \textbf{$V_{sum}$} denotes the accumulated magnitude of SLO violation, calculated as the sum of $max(c^d - c^*, 0)$ for all the minutes during the experiment period. It is important to note that a high $S_{vr}$ may be accompanied by a low $V_{sum}$; therefore, lower metric values signify desirable SLO satisfaction;
    \item \textbf{$R_{avg}$} denotes the cost of resources, calculated as the average number of installed nodes $x^d$ during the experiment period;
\end{itemize}

\subsubsection{Experimental Results for Question 2} \label{sec: Q2}

\begin{table*}[tb]
\caption{Comparison on four types of autoscalers. Each autoscaler is tested with sets $S_A$, $S_B$, and $S_C$. 
Three metrics are evaluated under different parameter settings.  
Lower metric values are desired, and the best are bolded.}
\begin{center}
\begin{threeparttable}
\begin{tabular}{llllllllllll}
    \cmidrule(r){1-12}
    \multirow{2}{*}{\textbf{Type}} & \multirow{2}{*}{\textbf{Name}} & \multirow{2}{*}{\textbf{Param}}
    & \multicolumn{3}{c}{\textbf{$S_A$}\tnote{a}} & \multicolumn{3}{c}{\textbf{$S_B$}\tnote{a}} & \multicolumn{3}{c}{\textbf{$S_C$}\tnote{a}}  \\
    \cmidrule(r){4-6} \cmidrule(r){7-9} \cmidrule(r){10-12}
     & & & $S_{vr}(\%)$ &\textbf{$V_{sum}$} &\textbf{$R_{avg}$} &$S_{vr}(\%)$ & \textbf{$V_{sum}$} &\textbf{$R_{avg}$} &$S_{vr}(\%)$ & \textbf{$V_{sum}$} &\textbf{$R_{avg}$} \\
    \cmidrule(r){1-3} \cmidrule(r){4-6} \cmidrule(r){7-9} \cmidrule(r){10-12}
    
    Reactive & Autopilot
    & $S = 90\%$ & 6.1 & 1.87 & 115.9  & 13.8 & 45.85 & \textbf{193.0}    & 9.6 & 34.52 & 132.6 \\
    && $S = 95\%$ & 2.8 & 0.87 & 119.0  & 1.8 & 4.98 &  266.7  & 4.3 & 9.64 & 158.8 \\
    \cmidrule(r){1-3} \cmidrule(r){4-6} \cmidrule(r){7-9} \cmidrule(r){10-12}
    
    Proactive & Madu
    &  $D=1$ & 13.7	  &  8.87 & \textbf{99.0}  & 5.6 & 13.20 & 214.9  & 7.0 & 8.74 & \textbf{69.5}   \\
    &&  $D=11$ & 13.6 &  8.65  & 99.1   & 2.7   & 3.22 & 252.5  &  6.5 & 8.34 & 70.0  \\
    \cmidrule(r){1-3} \cmidrule(r){4-6} \cmidrule(r){7-9} \cmidrule(r){10-12}
    
    Hybrid & HAS 
    &  $D=1$ & 1.6 & 0.73 &  105.1   & 5.3 & 13.96 & 223.0 & 3.6 & 3.36 & 81.3   \\
    &&  $D=11$ & 1.5 & 0.58 &104.6  & 2.2 & 2.21 & 263.2 & 2.8 & 3.13 & 80.4   \\
    \cmidrule(r){1-3} \cmidrule(r){4-6} \cmidrule(r){7-9} \cmidrule(r){10-12}
    
    Collaborative & \textbf{OptScaler} &  $D=1$ & 1.2	& 0.26  & 113.0  & 2.1 & 3.56 & 226.3  &1.1 & 0.63 & 84.8    \\
    &  \textbf{(proposed)} & $D=11$ & \textbf{0.7} & \textbf{0.25} & 112.3  & \textbf{0.6} & \textbf{0.61} & 265.7 &\textbf{0.7} & \textbf{0.55} & 90.5   	\\
    \cmidrule(r){1-12}
    
    \end{tabular}
    \begin{tablenotes}
        \footnotesize
        \item[a] Due to the randomness in sampling the CPU utilization feedback, all evaluation metrics are averaged over five repetitive runs and calculated in minute-level resolution. $S_{vr}$ denotes the SLO violation rate, $V_{sum}$ denotes the accumulated magnitude of SLO violation, and $R_{avg}$ denotes the resource cost.
    \end{tablenotes}
\end{threeparttable}
\label{tab:azure_offline_result}
\end{center}
\end{table*}

The statistical results are summarized in \autoref{tab:azure_offline_result}, and the metrics are explained in Section~\ref{sec: Metrics}. From the table, OptScaler excels in SLO satisfaction while maintaining relatively low resource cost, particularly for complex workloads (e.g., $S_B$ and $S_C$). Further details are provided below:
\begin{itemize}
    \item Under the quantile of $S =90\%$, Autopilot demonstrates an unacceptable level of safety across all sets; with $S =95\%$, it achieves an averaged performance in preventing SLO violations across all experiment groups, but the resource cost $R_{avg}$ is notably higher than in other experiment groups, particularly in $S_A$ and $S_C$, potentially limiting its applicability under tight resource budgets. The high magnitude of $V_{sum}$ presents another issue for Autopilot to address;
    
    \item Madu demonstrates high risks in terms of both SLO violation rate $S_{vr}$ and magnitude $V_{sum}$. This can be attributed to the disparity between estimated CPU utilization and the real system feedback. Without adjustment by a reactive module, this gap will consistently undermine the reliability of final scaling decisions. This underscores the assertion in \cite{straesser2022not} that a reactive scaling component should be a part of every production-ready autoscaler, complementing an unsatisfactory proactive counterpart;

    \item OptScaler outperforms Autopilot in terms of both safety and resource costs. When compared to the hybrid autoscaler HAS with $D=11$, OptScaler incurs a slightly higher resource cost $R_{avg}$ (1\% $\sim$ 12\% more) to achieve a significantly lower SLO violation rate $S_{vr}$ (53\% $\sim$ 75\% less) and magnitude $V_{sum}$ (56\% $\sim$ 82\% less). This advantage is even more pronounced when compared to Autopilot and Madu. It demonstrates OptScaler's superiority in SLO satisfaction, particularly for challenging workload patterns such as in sets $S_B$ and $S_C$;

    \item The comparison between $D = 1$ and $D = 11$ reveals that the multi-timestep optimization structure in MPC brings a consistent improvement of the rate $S_{vr}$ and magnitude $V_{sum}$ of SLO violation in both OptScaler and HAS.
\end{itemize}

\begin{figure*}[t]
	\centering
	\includegraphics[width = 2\columnwidth]{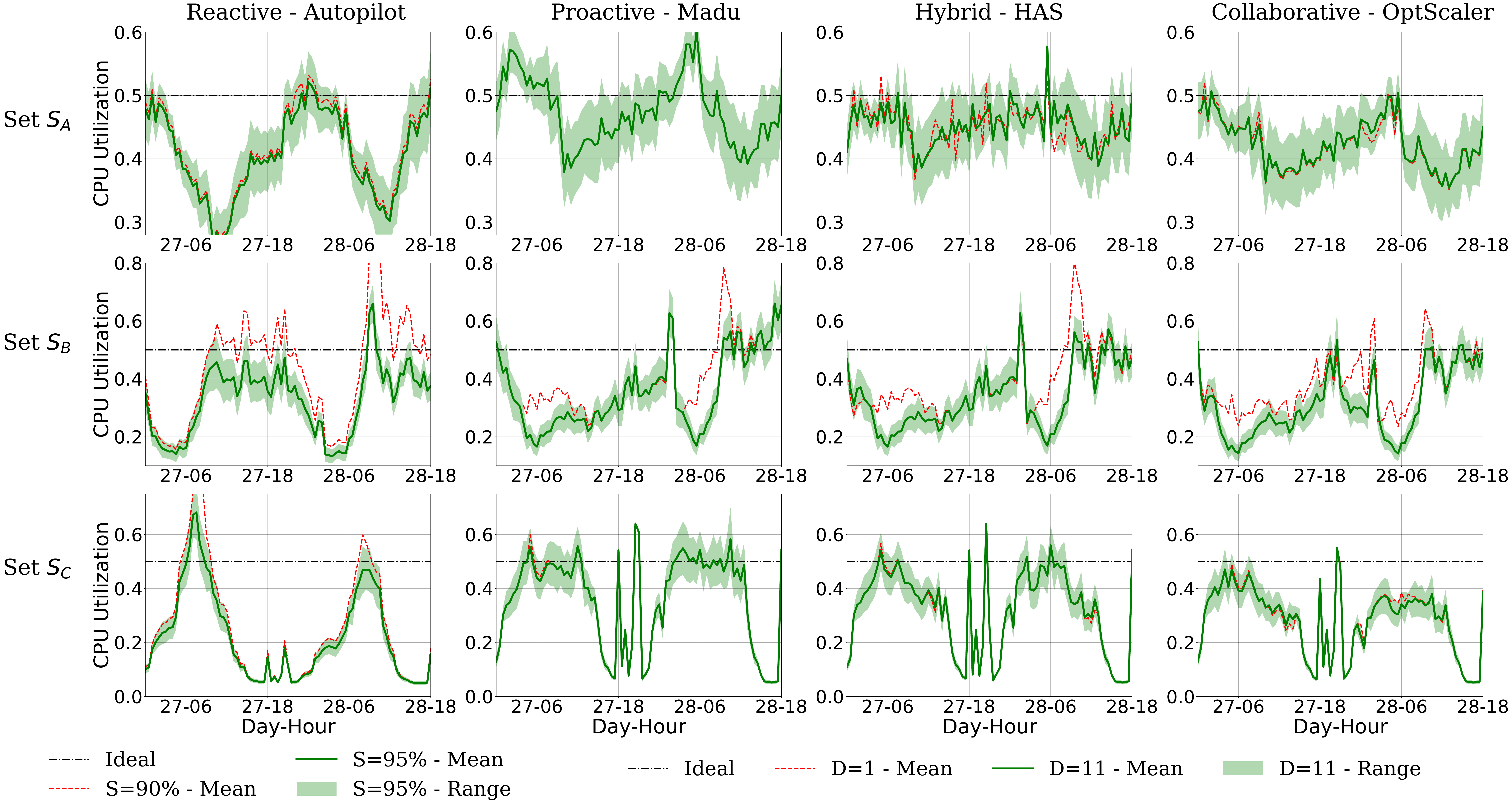}
	\caption{Experimental results of CPU utilization of four types of autoscalers.
    Each autoscaler (column) is tested with sets $S_A$ (top row), $S_B$ (middle row), and $S_C$ (bottom row) under two different parameter settings (as indicated by the bottom legends). 
    }
    \Description{}
    \label{fig:set_CPU_result}
\end{figure*}

We have further generated \autoref{fig:set_CPU_result} to present the experimental details of CPU utilization for clusters linked with sets $S_A$, $S_B$, and $S_C$. In this figure, we display the \textbf{Ideal} line for $c^*$, along with two \textbf{Mean} lines representing actual CPU utilization under different parameter settings. The light \textbf{Range} area illustrates the $\alpha=0.95$ confidence interval of the actual CPU utilization. For improved clarity, we only depict the Range area for the line demonstrating better performance, as indicated in \autoref{tab:azure_offline_result}; specifically, we display $S =95\%$ for Autopilot and $D = 11$ for the three optimization-based autoscalers. Any part of the line above the Ideal line signifies the risk of SLO violations. From \autoref{fig:set_CPU_result}, we observe four key insights:

\begin{itemize}
    \item Autopilot struggles with significant fluctuations in CPU utilization due to the absence of workload prediction. As showcased for $S_C$ in \autoref{fig:Azure_series_new} and \autoref{fig:set_CPU_result}, from 6 AM to 9 AM on Jul 27 (indicated as 27-06 to 27-09 on the x-axis), Autopilot exhibits delayed responsiveness to the substantial surge in workloads, potentially leading to severe SLO violations. Using $S =95\%$ only serves as a partial solution. The underlying problem is that Autopilot only considers a limited historical window in preparing resources for the future. Consequently, when a more significant workload spike occurs, the risk of violations could reoccur; 
    
    \item The pure proactive scaler Madu also experiences SLO violations due to the absence of a reactive module to rectify scaling errors, leading to a wider range of SLO violations than HAS and OptScaler, as evidenced in the period from 28-06 to 28-12 for $S_C$ in \autoref{fig:set_CPU_result}. This presents an undesirable outcome for cloud service providers. In contrast, the hybrid scaler HAS demonstrates less frequent SLO violations, underscoring the effectiveness of including both proactive and reactive modules in ensuring SLO compliance;
    
    \item OptScaler outperforms HAS by maintaining a safer fluctuation that significantly reduces both the duration and magnitude of SLO violations, as also reflected in $S_{vr}$ and $V_{sum}$ in \autoref{tab:azure_offline_result}. For example, for set $S_C$ in \autoref{fig:Azure_series_new} and \autoref{fig:set_CPU_result}, during two workload spikes around 6 PM and 9 PM (27-18 and 27-21 on the x-axis), HAS experiences two instances of SLO violations with CPU utilization reaching 0.63, while OptScaler encounters a single violation with less severe impact at 0.51;
    
    \item OptScaler, HAS, and Madu exhibit similar overall trends, but OptScaler consistently maintains CPU utilization at a much safer level, particularly under $D = 11$. Interestingly, the increase in the time horizon $D$ from $1$ to $11$ in MPC-based optimization results in less satisfactory CPU utilization for HAS and Madu compared to OptScaler. This is ascribed to their CPU estimator's inability to adapt to the evolving cloud environment characterized by less regular workload patterns. This illustrates the significance of incorporating a reactive module to align with MPC-based optimization, as implemented in OptScaler.
\end{itemize}

In summary, Autopilot exhibits limitations in effectively managing substantial workload fluctuations, owing to its hysteretic nature. However, when considering the three autoscalers equipped with a proactive module, OptScaler, as a collaborative autoscaler, demonstrates the most resilient performance at a slightly higher resource cost, followed by HAS, the hybrid autoscaler, with Madu, a pure proactive autoscaler, ranking at the bottom. This observation suggests that OptScaler effectively strikes a balance between resource cost and the risk of SLO violations, particularly for the co-located LRAs characterized by complex workload patterns.

\subsubsection{Experimental Results for Question 3} \label{sec: Q3}
We compare OptScaler with the existing HAS framework (the second-best performer in \autoref{tab:azure_offline_result}) using NBEATS (the second-best prediction model in \autoref{tab:ts_pred_res}). We design the experiments in an ablative manner, i.e., we adopt NBEATS to replace the workload prediction model in OptScaler to illustrate the impact of different prediction models on the final decision. The results are presented in \autoref{tab:compare_OptScaler}. 

Upon examining \autoref{tab:compare_OptScaler}, a comparison between Exp.2 and Exp.1 reveals a relatively $27\% \sim 60\%$ improvement of $S_{vr}$ by upgrading the autoscaling framework alone. Comparing Exp.3 to Exp.2 demonstrates a further improvement of $25\% \sim 53\%$ in $S_{vr}$, along with significant reductions in both $V_{sum}$ and $R_{avg}$, achieved by advancing the prediction model from NBEATS to our model. Consequently, when comparing OptScaler as a whole (Exp.3) to a best hybrid framework (HAS with NBEATS in Exp.1), OptScaler incurs $4\% \sim 15\%$ more resource costs in order to achieve $36\% \sim 81\%$ fewer SLO violations. This demonstrates the clear benefits of OptScaler for LRAs with demanding SLOs.
In essence, \autoref{tab:compare_OptScaler} serves to illustrate that the robustness of OptScaler is derived from enhancements in both the workload prediction model and the autoscaling framework.

\begin{table*}[tb]
\caption{Comparison between OptScaler and the existing hybrid framework of combining NBEATS and HAS. Ablative studies are conducted on all sets under $D = 11$. Metrics with lower values are desired, and the best are bolded.}

\begin{center}
\begin{tabular}{llllllllllll}
    \cmidrule(r){1-12}
    \multirow{2}{*}{\textbf{ID}} & \multirow{2}{*}{\makecell[l]{\textbf{Autoscaling} \\ \textbf{Framework}} } & \multirow{2}{*}{\makecell[l]{\textbf{Prediction} \\ \textbf{Model}}}
& \multicolumn{3}{c}{\textbf{$S_A$}\tnote{a}} & \multicolumn{3}{c}{\textbf{$S_B$}\tnote{a}} & \multicolumn{3}{c}{\textbf{$S_C$}\tnote{a}}  \\
    \cmidrule(r){4-6} \cmidrule(r){7-9} \cmidrule(r){10-12}
     & & & $S_{vr}(\%)$ &\textbf{$V_{sum}$} &\textbf{$R_{avg}$} &$S_{vr}(\%)$ & \textbf{$V_{sum}$} &\textbf{$R_{avg}$} &$S_{vr}(\%)$ & \textbf{$V_{sum}$} &\textbf{$R_{avg}$} \\
    \cmidrule(r){1-3} \cmidrule(r){4-6} \cmidrule(r){7-9} \cmidrule(r){10-12}

    Exp.1  & HAS & NBEATS & 1.1 & 0.60 & \textbf{105.9} & 2.0 & 3.04  & \textbf{254.3} & 3.7 & 8.75 & \textbf{78.8} \\
    \cmidrule(r){1-3} \cmidrule(r){4-6} \cmidrule(r){7-9} \cmidrule(r){10-12}
    
    Exp.2 & OptScaler & NBEATS & 0.8 & 0.40  & 114.5  & 0.8 & 0.76 & 277.6 & 1.5 & 5.32 & 133.8  \\
    \cmidrule(r){1-3} \cmidrule(r){4-6} \cmidrule(r){7-9} \cmidrule(r){10-12}
    
    Exp.3 & \textbf{OptScaler} & \textbf{Our model}  & \textbf{0.7} & \textbf{0.25} &  112.3 & \textbf{0.6}& \textbf{0.61} & 265.7 & \textbf{0.7} & \textbf{0.55} & 90.5  \\
    \cmidrule(r){1-12}
\end{tabular}
\label{tab:compare_OptScaler}
\end{center}
\end{table*}

\subsubsection{Further discussions} \label{sec: Q2_further}
The experimental results above provide a broad response to Question 2 and 3. However, in practice, the scenario may vary based on different scaling frequencies and unforeseen workloads. Therefore, in this subsection, we further enrich our investigation by delving into the following specific situations.

\begin{figure}[tb]
	\centering
	\includegraphics[width = 0.9\columnwidth]{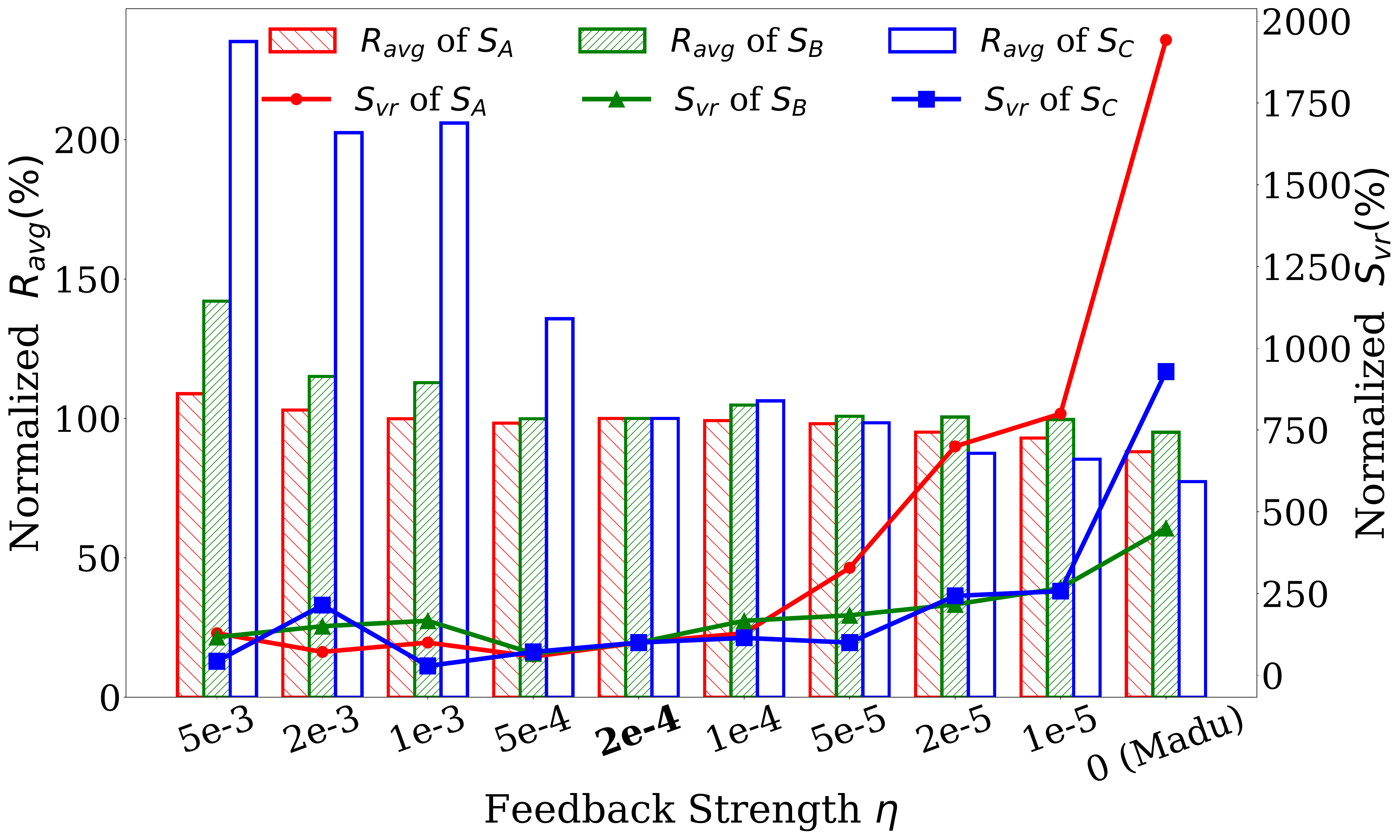}
	\caption{Sensitivity experiments of feedback strength $\eta$ in OptScaler show a clear balance between resource cost $R_{avg}$ (plotted as bars) and SLO violation rate $S_{vr}$ (plotted as lines). Metrics are evaluated under $D = 11$ and normalized with respect to the results of $\eta=2e-4$, the lower values the better.}
        \Description{}
	\label{fig:sensi_analysis}
\end{figure}

\paragraph{Sensitivity analysis} In order to explore the impact of hyperparameters in OptScaler, we conducted sensitivity experiments on a key parameter, $\eta$, which represents the strength of feedback in the reactive module  (refer to Algorithm \ref{alg:olr}). We examined various values of $\eta$, ranging from $5e-3$ to $0$ in increments of order of magnitude. When $\eta=0$, it effectively simulates the closed state of the reactive module (comparable to Madu). The results are presented in \autoref{fig:sensi_analysis}, where the left and right y-axes depict the crucial metrics of resource cost $R_{avg}$ (shown as bars) and SLO violation rate $S_{vr}$ (depicted as lines), respectively. Lower values are favourable for both metrics. We normalized each metric as a percentage of the corresponding result with $\eta=2e-4$ (the default value of OptScaler),  allowing for a comparison of performance across different $\eta$ values.

The findings portrayed in \autoref{fig:sensi_analysis} reveal a distinct trend wherein $R_{avg}$ and $S_{vr}$ move inversely as $\eta$ decreases from $5e-3$ to 0. Thus, it becomes evident that $\eta$ significantly influences the trade-off between the two evaluation metrics. Opting for a higher feedback strength with $\eta$ proves effective in suppressing SLO violations. However, an overly aggressive reactive strategy with a large $\eta$ also leads to an over-provisioning of resources. \autoref{fig:sensi_analysis} suggests that a relatively optimal range for $\eta$ falls between 5e-4 to 1e-4 across all sets. It is important to recognize that in real-world implementation, the selection of $\eta$ should be fine-tuned according to the local policies governing SLOs and resource budgets.

\begin{table}[tb]
\caption{Comparison on the response of OptScaler and Autopilot to unpredictable workloads on set $S_A$. Lower metric values are desired, and the best are bolded.}
\begin{center}
\begin{threeparttable}
\begin{tabular}{llllll}
    \cmidrule(r){1-6}
    \multirow{2}{*}{\textbf{Name}} & \multirow{2}{*}{\textbf{Param}}
    & \multicolumn{2}{c}{\textbf{30-minute spike}} & \multicolumn{2}{c}{\textbf{3-hour spike}}  \\
    \cmidrule(r){3-4} \cmidrule(r){5-6}
     & & $S_{vr}(\%)$  &\textbf{$R_{avg}$} &$S_{vr}(\%)$ &\textbf{$R_{avg}$} \\
    \cmidrule(r){1-2} \cmidrule(r){3-4} \cmidrule(r){5-6}
    
    Autopilot
    & $S = 95\%$ & 3.7 &	\textbf{120.9} &  12.7 & 	567.0  \\
    & $S = 99\%$ & 2.3 &	569.6 &  12.6	& 609.7 \\
    \cmidrule(r){1-2} \cmidrule(r){3-4} \cmidrule(r){5-6}
    
    \textbf{OptScaler}   
    &  $D=1$  & 2.4  &  306.4 & 12.6 &  \textbf{316.7}      \\
    &  $D=11$ & \textbf{2.2} &  308.3 & \textbf{12.5} &  317.6      	\\
    \cmidrule(r){1-6}
    
    \end{tabular}
\end{threeparttable}
\label{tab:azure_spike}
\end{center}
\end{table}

\paragraph{Unpredictable workload} To further evaluate the robustness of OptScaler, we mock some totally unpredictable spikes by directly increasing the workloads to 10 times of the original value based on $S_A$, and test the response of OptScaler and its reactive competitor Autopilot to that workload. Note that Autopilot's ability to recognize the unprecedented spike depends on a percentile parameter $S$ as mentioned in Section \ref{sec: Autoscaling_Frameworks}. A higher $S$ makes it more sensitive to unpredictable workload. Our results in \autoref{tab:azure_spike} show that when $S=95\%$, Autopilot tends to ignore the unprecedented spike, especially when the duration of the spike is short (e.g. 30-minute). In contrast, OptScaler recognizes the event and scales up the resources accordingly, resulting in a lower SLO violation $S_{vr}$; when $S$ continues to increase (e.g. 99\%), both OptScaler and Autopilot can quickly respond to the event; however, a higher $S$ will lead to an over-conservative Autopilot, and OptScaler leads to significantly less resource cost $R_{avg}$ than Autopilot does when they achieve a similar level of $S_{vr}$.

\paragraph{Scaling frequency} In practice, clusters may vary in scaling frequency, hence we enrich our experiments by investigating autoscaling frameworks under different intervals $h$. Here we test $h = 10$ for 10-minute interval, and $h = 60$ for 1-hour interval. We omit tables for lack of space, and the results are summarized as follows: 1) OptScaler achieves the best $S_{vr}$ and $V_{sum}$ under both $h$, followed by HAS as the second-best, which aligns with that of $h=30$; 2) Under the same horizon $D$, OptScaler generally costs more resources in the scenario of 1-hour interval than that of 10-minute, as it considers a longer period of future workloads; 3) The evaluation metrics of Autopilot shows no significant differences under different $h$, as it has no proactive or optimization module, and $h$ has very limited impact on the evolving process of scaling metrics.

\section{Deployment} \label{sec: Deploy}

\begin{figure}[tb]
	\centering
	\includegraphics[width = 0.75\columnwidth]{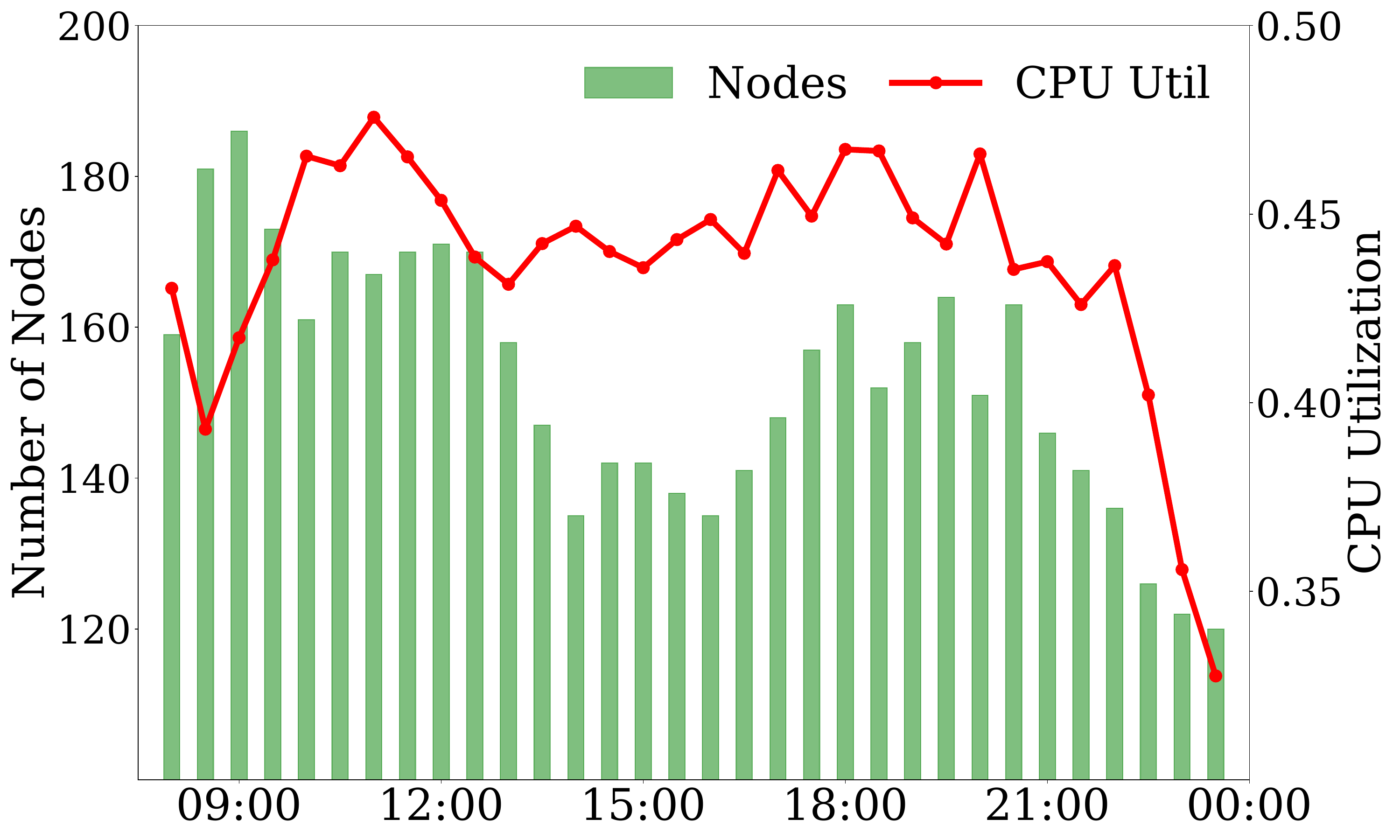}
	\caption{Online experimental results of OptScaler on the number of installed nodes (bars with left y-axis) and CPU utilization (red line with right y-axis). OptScaler demonstrates its ability to regulate CPU utilization.}
        \Description{}
	\label{fig:online_res}
\end{figure}

Before deployment, we tested OptScaler in a production cluster with co-located LRAs to verify its effectiveness in the online environment. We refrained from comparing OptScaler with other frameworks due to concerns about their SLO violations exceeding 2\%, as indicated by $S_{vr}$ in \autoref{tab:azure_offline_result}. We conducted our experiment from 8 AM to 12 PM. From the result shown in \autoref{fig:online_res}, OptScaler achieved a steady CPU utilization (between 0.4 and 0.48) during the whole daytime. It started to decrease after 10 PM as the total workload dropped to half of that in the daytime. Besides, the previous resource cost was fixed at $190$ to ensure safety, and the cost $R_{avg}$ by OptScaler was $152.9$, saving 19.5\% of cloud resources.

In addition, OptScaler benefits users for its interpretability. It allows users to trace back to specific conditions and parameters in the model to understand why a scaling decision is made. When an unexpected scaling occurs, this model-based system can facilitate the debugging and make further improvement. Hence, users are likely to trust OptScaler more than other black-box scaling methods. 

\balance

OptScaler has been successfully deployed as an integrated platform of prediction and decision-making at Alipay, supporting the autoscaling of more than 100 online LRAs. To access the platform of OptScaler, a cluster typically undergoes the following steps: 1) The user configures necessary parameters for the cluster and the algorithm; 2) The cloud system monitors and collects the workloads and metrics of the cluster, and the data will be saved to a database; 3) OptScaler loads historical workload data for the proactive module to train model offline at a regular pace, and outputs the prediction results at real-time; 4) OptScaler loads the latest system feedback from the database to fine-tune the reactive module; 5) The optimization module makes a final scaling decision, which is returned to the cloud cluster for implementation.

\section{Conclusion} \label{sec: Conclude}
We present OptScaler as a pioneering collaborative autoscaling framework, designed as an upgrade to the widely used hybrid framework. OptScaler stands out as the first framework to address the collaboration of proactive and reactive methods through sophisticated optimization techniques, such as Model Predictive Control (MPC). Unlike approaches that employ proactive and reactive methods independently, which can lead to incompatibilities, OptScaler orchestrates the strengths of both proactive and reactive modules to effectively manage workload fluctuations and system uncertainty, significantly enhancing the robustness of cloud clusters. OptScaler surpasses other prevalent autoscaling frameworks thanks to its superior prediction model and collaborative mechanism. Offline experiments demonstrate that OptScaler effectively mitigates the risk of SLO violations by a minimum of 36\% in comparison to other frameworks. In online experiments, OptScaler maintains desirable CPU utilization while achieving substantial savings of up to 19.5\% in cloud resources. It is noteworthy that OptScaler has already been deployed online to support the autoscaling of LRAs at Alipay, a leading global payment platform.

\bibliographystyle{ACM-Reference-Format}
\bibliography{VLDB24}

\end{document}